\documentclass{jors}

\definecolor{mygray}{gray}{0.6}

\usepackage{longtable}
\begin{document}

{\bf Software paper for submission to the Journal of Open Research Software} \\

To complete this template, please replace the blue text with your own. The paper has three main sections: (1) Overview; (2) Availability; (3) Reuse potential. \\

Please submit the completed paper to: editor.jors@ubiquitypress.com

\rule{\textwidth}{1pt}

\section*{(1) Overview}

\vspace{0.5cm}

\section*{ArcLP: A Matlab implementation of an $\mathcal{O}(\sqrt{n}L)$ 
arc-search infeasible interior-point algorithm for linear programming}


\section*{Yang, Yaguang;}


\section*{Yaguang Yang is an independent researcher and the sole author of the paper.}

\section*{Abstract}

This paper presents a Matlab implementation of an arc-search infeasible interior-point algorithm for linear programming (LP), which has a proven polynomial bound of $\mathcal{O}(\sqrt{n}L)$, the best among all interior-point algorithms for LP. Software architecture and major functions are discussed. Its ease of use is described by a simple example. Crucial strategies are summarized. Quality of the software is assured because this software has been extensively tested on both PC and Linux for the widely used Netlib benchmark linear programming problems in standard form. Some benchmark test problems involve tens of thousands of constraints and hundreds of thousands of variables.  For all tested problems, the code found the optimal solution. The numerical results have been compared to those obtained by the popular Mehrotra’s predictor-corrector algorithm.
We conclude that the implemented algorithm not only has the best polynomial bound but also is computationally competitive compared to the popular Mehrotra’s predictor-corrector algorithm.

\section*{Keywords}
linear programming; infeasible interior-point algorithm; arc-search

\section*{Introduction}
Linear programming is one of the most studied 
mathematical problems. It has been used in 
almost every scientific discipline. There are two
popular methods to compute the solution of the 
linear programming problem: one is the famous 
simplex method developed in 1949 by
Dantzig in \cite{dantzig49}; the other one is 
the interior-point method that became popular 
after Karmarkar showed that his interior-point
algorithm converges in polynomial time \cite{karmarkar84}. 
Since simplex-method-based algorithms 
do not converge in polynomial time in the worst 
case \cite{psz08}, theoretically, interior-point method
is more attractive because many interior-point algorithms based on
interior-point method converge in polynomial time
in the worst case. In addition, Mehrotra 
\cite{Mehrotra92} developed his predictor-corrector interior-point algorithm, 
which has been demonstrated to be competitive in numerical tests to the efficient 
simplex method for large-scale problems \cite{lms92}. However, unlike other 
interior-point algorithms which are less efficient in numerical test than
Mehrotra's algorithm, the convergence of Mehrotra's predictor-corrector algorithm has not been proved. As a matter of fact, examples were found that
Mehrotra's algorithm does not converge, therefore, 
efficient interior-point algorithms that converge in 
polynomial time were proposed in \cite{yang18}.

In this paper, we present an implementation of an algorithm in 
\cite{yang18}, which has a proven best polynomial bound 
$\mathcal{O}(\sqrt{n}L)$ among all interior-point 
algorithms developed for linear programming problems. 
Yet numerical test results demonstrate that this algorithm 
is more efficient and robust than Mehrotra's algorithm
which was regarded as the most efficient interior-point 
algorithm for
linear programming problem for many years \cite{todd02,wright95}.

Mathematically, the linear programming in standard form is expressed 
as follows:
\begin{eqnarray}
\min \hspace{0.05in} c^{T}x, \hspace{0.15in} \mbox{\rm subject to} 
\hspace{0.1in}  Ax=b, \hspace{0.1in} x \ge 0,
\label{LP}
\end{eqnarray}
where $A \in R^{m \times n}$, $b \in R^{m} $, $c \in R^{n}$ 
are given, and $x \in R^n$  is the vector to be optimized. 
The dual of (\ref{LP}) is:
\begin{eqnarray}
& \max \hspace{0.05in} b^{T}\lambda, \hspace{0.15in} 
\mbox{\rm subject to} 
\hspace{0.1in}  A^{T}\lambda+s=c, \hspace{0.1in} s \ge 0,
\label{DP}
\end{eqnarray}
with dual variable vector $\lambda \in R^{m}$, and 
dual slack vector $s \in R^{n}$. To solve this problem, 
interior-point algorithms start from an interior-point of 
the problem which meets the condition of $(x^0,s^0)>0$; 
all iterates are interior-points, i.e., $(x^k,s^k)>0$, 
where the superscript $k$ stands for the $k$th iteration; 
the search is to find the solution that meets the KKT conditions,
the necessary and sufficient conditions for the optimal
solution of (\ref{LP}).

There are many commercial software packages such 
as {\tt OB1} \cite{lms92}, 
{\tt CPLEX} \cite{nssb22}, and {\tt LIPSOL} \cite{mathworks,zhang96} that can be used to solve linear programming problem
(\ref{LP}). All these software packages that use IPM for LP are based 
on Mehrotra's algorithm \cite{wright95} which, as we have pointed out, 
does not have a convergence result. In contrast, 
the software presented in this paper provides not only
a free code that solves the same problem, the 
code is also based on a recent algorithm that is proved
to converge in polynomial time with the best polynomial bound. 
Given the fact that many application problems 
in science and engineering can be formulated as or reduced to a 
linear programming problem, this software provides 
a superior alternative for 
potential users with broad backgrounds and interests.

All the test problems are from Netlib benchmark 
set in MAtlab format that was widely used by people in optimization 
community to test and evaluate the code designed
for linear programming problems \cite{NETLIB}. 
We demonstrated in \cite{yang18} the computational merits of the 
implemented code {\tt arcLP} by testing it 
along with a Matlab code that implements Mehrotra’s 
algorithm using Netlib problems in Matlab format and comparing 
the test results. To have a fair comparison, we used 
the same initial point, the same pre-process and 
post-process, and the same termination criteria 
for the two algorithms for all test problems. This result is provided in \cite{yang18}.

Given problem data $A, b, c$, and an option 
$d$ ($d=1$ uses a function to deal with degenerate 
problems which will be discussed in the next
section, and $d=0$ is the default option which will not 
use the function),
we use 
\newline\hspace*{\fill}
[x,obj,kk,infe,lambda,s,exflag]=arcLP(A,b,c,d,tol,iter)
\hspace*{\fill}\newline
to call the implemented code {\tt arcLP}.
The output $x$ is the optimal solution, $obj$ is the optimal objective value, $kk$ is the number of iterations used to find the optimal solution, $infe=\| Ax-b \|_2$ is the residual of the equality constraints, $lambda$ is the dual optimal solution, $s$ is the corresponding value of the slack variable, and $exflag$ is the exit flag returned by the function, with $exflag=0$ for success termination, $exflag=1$ for a certificate of an infeasible instance, $exflag=2$ for a certificate of an unbounded instance, $exflag=3$ for a certificate that both the primal and dual problems are infeasible, and $exflag=4$ for an incomplete input data $(A,b,c)$.

\section*{Implementation and architecture}
\label{Software-Architecture}

The software has several functions: 

\begin{itemize}
\item [] function [x,obj,kk,infe,lambda,s,exflag]=arcLP(A,b,c,d,tol,iter)
is the main function, which includes pre-process, an 
option to make matrix $A$ full rank if it is not, 
intuitive selection of initial point, the main algorithm, 
the post-process, and a final check of optimality. 

\item [] function [minA,maxA]=calRatioCondition(A)
finds the smallest and biggest absolute 
non-zeros of $A$. The information is used to decide 
if we need to scale the matrix $A$ , vector $b$, and vector $c$.

\item [] function 
makeAfull(A,b,c,x,s,lambda,rB,rC,lOld,xOld,sOld) uses 
Markowitz's pivot criterion \cite{dobes05} to remove 
dependent rows of $A$ to make the rows of matrix 
$A$ linearly independent while not losing the sparsity.

\end{itemize}

%
Some details and strategies that are used to enhance the efficiency and robustness are described below.
\begin{itemize}
\item{\bf Initial point selection}

Initial-point selection has long been recognized as an important factor affecting the computational efficiency of most infeasible interior-point algorithms \cite{cmww97, zhang96}. In this work, we employ the two methods proposed in \cite{Mehrotra92,lms92} to generate candidate initial points. We then evaluate the quantity
\begin{equation}
\max \{ \| A x^0-b\|, \| A^{T} \lambda^0 +s^0 -c \|, x^{0^{T}}s^0/n \}
\label{initpoint}
\end{equation}
for each candidate and select the initial point that yields the smaller value. This choice is motivated by the expectation that a smaller value of the above measure may lead to a reduction in the number of iterations required by the algorithm (see \cite{yang16} for further details).

\item{\bf Pre-process and Post-process}\label{presolver}

Preprocessing strategies for linear programming problems in the form of (\ref{LP}) were thoroughly investigated in \cite{yang16}. These strategies enable the solver to handle problem instances whose data do not satisfy the assumptions imposed in \cite{yang18}. In the present paper, we employ the same set of preprocessing procedures as those developed in \cite{yang16}. The postprocessing procedures are likewise identical to those described in \cite{yang16}.

\item{\bf Matrix scaling}

Matrix scaling was originally developed to mitigate the effects of ill-conditioned matrices. However, based on the tests and analysis reported in \cite{yang16}, matrix scaling was found not to improve computational efficiency in general. Therefore, scaling is not included in our implementation. Nevertheless, the ratio
\begin{equation}
\frac{\max{|A_{i,j}|}}{\min { \{ | A_{k,l} |  A_{k,l} \ne 0 \} } }
\label{ratio}
\end{equation}
is used to determine whether one of the preprocessing rules proposed in \cite{yang16} should be applied.

\item{\bf Removing row dependency from $A$}

The removal of row dependencies in $A$ was studied in \cite{andersen95}, where Andersen proposed an efficient method for detecting and eliminating dependent rows. However, based on the study presented in \cite{yang16}, we have chosen not to incorporate this procedure into our implementation. An exception may be made if it becomes necessary as part of the strategy for handling degenerate solutions, as discussed below.

\item{\bf Linear algebra for sparse Cholesky factorization}

Similar to Mehrotra's algorithm, the dominant computational cost of the proposed algorithms lies in solving sparse Cholesky systems, which can be expressed in the abstract form
\begin{equation}
AD^2A^{T} u = L \Lambda L^{T} u = v,
\label{useLater}
\end{equation}
where $D$ and $\Lambda$ are diagonal matrices, $L$ is a lower triangular matrix, and $u$ and $v$ are vectors. Many widely used LP solvers \cite{cmww97,zhang96} rely on the software package \cite{np93}, which incorporates linear algebra techniques specifically designed for sparse Cholesky factorization \cite{liu85}. However, MATLAB does not provide all of the features required to robustly handle ill-conditioned matrices in this context. Therefore, we adopt the implementation strategy described in \cite{yang16}.

\item{\bf Handling degenerate solutions}

The difficulties caused by degenerate solutions in interior-point algorithms for linear programming have long been recognized \cite{ghrt93}. Similar observations were reported in \cite{gmstw86}. To address this issue, our implementation includes an optional procedure for handling degenerate solutions, following the approach described in \cite{yang16}.

\item{\bf Analytic solution for step angle of $\alpha_k$}

The proposed arc-search algorithm employs a distinctive search strategy that follows an arc toward the optimizer rather than a straight line. The corresponding step angle can be computed analytically, and the derivation of this formula is presented in \cite{yang09,yang16}.

\item{\bf Select centering parameter $\sigma_k$}

The proposed interior-point algorithm differs from most existing interior-point methods in its selection of the centering parameter. Specifically, the parameter is determined optimally through a specialized Golden Section search that exploits a particular property of the problem.

It is well known that the classical Golden Section search reduces the interval length to approximately (0.618) times its previous value at each iteration \cite{luenberger84}. In contrast, the proposed algorithm reduces the interval length to (0.5) times its previous value in every iteration. Consequently, the proposed approach achieves a faster reduction of the interval and is therefore more efficient. The implementation used in this work follows the procedure described in \cite{yang18}.

\item{\bf Rescale step angle $\alpha_k$}

To enhance numerical stability, the step angle $\alpha_k$ obtained from the procedure described above is rescaled at each iteration according to $$\alpha_k = \min \{ 0.9999 \alpha_k, 0.99 \pi/2 \} < 0.99 \pi/2.$$ Our numerical experiments indicate that this rescaling prevents $x^k$ and $s^k$ from approaching zero too rapidly during the early iterations. As a result, it mitigates the numerical difficulties associated with solving (\ref{useLater}) and improves the overall robustness of the algorithm.

\item{\bf Detection infeasible and unbounded instances }

Some LP instances may be unbounded. In the interior-point method, unboundedness can be detected by checking whether $\max_i |x_i| > M, $ where $M = 10^{10}$ is the threshold implemented in {\tt arcLP}. Detecting infeasibility is somewhat more challenging. One approach is based on the duality theorem \cite[Theorem 13.1]{nocedal99}, if a dual variable $\lambda_i$ becomes unbounded, then the corresponding LP instance is infeasible. Consequently, infeasibility may be detected by checking whether $\max_i |\lambda_i| > M.$ Another indictor of infeasibility is when the algorithm is unable find a step size that maintains positivity of primal variable $x$ and slack variable $s$ while simultaneously reducing the duality gap. There are also cases in which both the primal and dual problems are infeasible. In such situations, {\tt arcLP} reports this status when the number of iterations reaches the prescribed maximum.

In fact, {\tt arcLP} may detect infeasible or unbounded instances during the presolve phase described in Section~\ref{presolver} (see \cite{yang16}). For example, if a column of $A$ consists entirely of zeros while the corresponding objective coefficient satisfies $c_i < 0$, then the LP is unbounded. In this case, the variable $x_i$ can increase without bound without affecting any of the constraints, causing the objective value to decrease without bound. When such infeasible or unbounded instances are identified during presolve, {\tt arcLP} reports the detected status and terminates without proceeding to the optimization phase.

We performed some numerical tests on these implemented heuristics and the results demonstrated that they provide reasonable detection of infeasible and unbounded instances.

\item{\bf Termination criteria}

The primary stopping criterion used in the implementation follows the standard convention adopted by most infeasible interior-point software packages, including LIPSOL \cite{zhang96}.
\[
\frac{\|r_b^k\|}{\max \lbrace 1, \| b\| \rbrace }
+\frac{\|r_c^k \|}{\max \lbrace 1, \| c\|  \rbrace }
+\frac{ \mu_k }{\max \lbrace 1, \| c^{T}x^k \|, \|b^{T}\lambda^k \|  \rbrace } 
< 10^{-8},
\] 
where $r_b^k=A x^k -b$ and $r_c^k=A^{T} \lambda^k + s^k -c$.

\item{\bf Miscellaneous}
In case that the algorithms fail to find a good search direction,
the program also stops if step sizes $\alpha_k^x < 10^{-8}$ and
$\alpha_k^s < 10^{-8}$. 

Finally, if (a) due to the numerical problem, $r_b^{k}$ or $r_c^{k}$ does not decrease but 
$10r_b^{k-1}<r_b^{k}$ or $10r_c^{k-1}<r_c^{k}$, or (b) if $\mu<10^{-8}$, the program stops.

\item{\bf A simple illustrative example}
\label{}

Let us consider 
\[
\min x_1, \hspace{0.15in} s.t. \hspace{0.1in} x_1+x_2 = 5,
\hspace{0.1in} x_1 \ge 0, \hspace{0.1in} x_2 \ge 0. 
\]
Therefore $A=[1 \hspace{0.08in} 1]$, $b=5$, and 
$c=[1 \hspace{0.08in} 0]^{T}$. 
Calling {\tt arcLP} gives the optimal solution 
$x^*=[0 \hspace{0.08in} 5]^{T}$
in five iterations.
\end{itemize}


\section*{Quality control}

We have tested the software extensively on both Linux and Windows computers for all Netlib benchmark linear programming problems in Matlab format represented in standard form.  The comprehensive test result is presented in \cite{yang18} which is very impressive (the largest number of constraints are 16675 and the largest number of variables are 104374 \cite{yang18}). 

Besides {\tt arcLP}, we also implemented the famous Mehrotra's
predictor-corrector algorithm as {\tt mehrotra}. We have tested the two
codes against all Netlib benchmark linear programming problems 
represented in standard form. These two codes adopted
the same pre-process and post-process, used the same 
initial points and options, and stopped with the
same criteria. The result is summarized in Table~\ref{tableIteration}.

\begin{center}
\begin{longtable}{|c|c|c|c|c|}
\caption{Test results of arclp.m and mehrotra.m for problems in Netlib}
\label{tableIteration}
\\   \hline    
Problem  & algorithm   & iter   &  obj    & infeasibility   \\
\hline

Adlittle
         &  mehrotra.m   & 15 &    2.2549e+05   & 3.4e-08    \\
         &    arcLP.m   & 16  &    2.2549e+05   &  3.0e-11  \\   \hline
Afiro
         &  mehrotra.m   &  9 & -464.7531   & 8.0e-12   \\
         &    arcLP.m   &  9 &  -464.7531   & 6.2e-13     \\  \hline
Agg
         &  mehrotra.m   &  22 & -3.5992e+07 & 5.2e-05  \\
         &    arcLP.m   &  20   &  -3.5992e+07   &   3.7e-06   \\ \hline
Agg2
         &  mehrotra.m   & 20  & -2.0239e+07 & 5.2e-07  \\
         &    arcLP.m   &  21 &   -2.0239e+07   &  3.1e-08    \\ \hline
Agg3 
         &  mehrotra.m   &  18 & 1.0312e+07  & 8.8e-09  \\
         &    arcLP.m   &  20 &   1.0312e+07  &   1.5e-08    \\ \hline
Bandm
         &  mehrotra.m   &  22 &  -158.6280   & 8.3e-10  \\
         &    arcLP.m   &  20 &  -158.6280  &  3.6e-11     \\ \hline
Beaconfd
         &  mehrotra.m   &  11 &  3.3592e+04  & 1.4e-10  \\
         &    arcLP.m   &  11 &  3.3592e+04  & 1.8e-12    \\ \hline
Blend  
         &  mehrotra.m   &   14 & -30.8122    & 4.9e-11  \\
         &    arcLP.m   &   14 &  -30.8122   &   1.6e-12    \\ \hline
Bnl1   
         &  mehrotra.m   &  35 &   1.9776e+03  & 3.4e-09  \\
         &    arcLP.m   & 34  &   1.9776e+03  &  2.9e-09   \\ \hline
Bnl2+ 
         &  mehrotra.m   &  38 &   1.8112e+03  & 9.3e-07  \\
         &    arcLP.m   &  35 &   1.8112e+03  &  3.5e-06   \\ \hline
Brandy 
         &  mehrotra.m   &  19  &  1.5185e+03  & 6.2e-08  \\
         &    arcLP.m   &  24  &  1.5185e+03  &  2.4e-06  \\ \hline
Degen2+
         &  mehrotra.m   &   17 &  -1.4352e+03  & 2.0e-10  \\
         &    arcLP.m   &  19  &  -1.4352e+03  &  5.9e-10  \\ \hline
Degen3*
         &  mehrotra.m   &   22 & -9.8729e+02   & 1.2e-09  \\
         &    arcLP.m   &  35  & -9.8729e+02   &  8.6e-08   \\ \hline
fffff800
         &  mehrotra.m   &   31 & 5.5568e+05  & 7.7e-04  \\
         &    arcLP.m   & 28   &    5.5568e+05   &   3.7e-09  \\ \hline
Israel
         &  mehrotra.m   &   29 & -8.9665e+05 & 1.8e-08  \\
         &    arcLP.m   &  27   & -8.9664e+05   &  3.4e-08  \\ \hline
Lotfi   
         &  mehrotra.m   &   18 & -25.2647    & 2.7e-07  \\
         &    arcLP.m   &  16  &  -25.2646 &  7.8e-09   \\  \hline
Maros\_r7
         &  mehrotra.m   & 21   &  1.4972e+06  & 6.4e-09  \\
         &    arcLP.m   &  20  &  1.4972e+06  &  1.7e-09  \\ \hline
Osa\_07+
         &  mehrotra.m   & 35 &   5.3578e+05  & 1.5e-07  \\
         &    arcLP.m   & 32 &   5.3578e+05  & 8.4e-10  \\ \hline
Osa\_14
         &  mehrotra.m   & 37 & 1.1065e+06    & 3.0e-08  \\
         &    arcLP.m   & 42 & 1.1065e+06    & 5.2e-09   \\ \hline
Osa\_30 
         &  mehrotra.m   & 36 &  2.1421e+06   & 1.3e-08  \\
         &    arcLP.m   & 42 &  2.1421e+06   &  1.3e-08    \\ \hline
Qap12 
         &  mehrotra.m   & 24 &  5.2289e+02   & 6.2e-09   \\
         &    arcLP.m   & 23 &  5.2289e+02    & 2.9e-10  \\ \hline
Qap15+ 
         &  mehrotra.m   & 44 &  1.0410e+03  & 1.5e-05   \\
         &    arcLP.m   & 28 &  1.0410e+03  &  8.4e-08  \\       \hline
Qap8+ 
         &  mehrotra.m   & 13 & 2.0350e+02  & 7.1e-09   \\
         &    arcLP.m   & 12 & 2.0350e+02  &  6.2e-11   \\ \hline
Sc105 
         &  mehrotra.m   &  11 & -52.2021    & 9.8e-11  \\
         &    arcLP.m   &  11 &    -52.2021   &   2.2e-12  \\ \hline
Sc205 
         &  mehrotra.m   &  12 & -52.2021    & 8.8e-11  \\
         &    arcLP.m   & 12  &   -52.2021     &   4.4e-11 \\     \hline
Sc50a 
         &  mehrotra.m   &  9  & -64.5751    & 8.3e-08  \\
         &    arcLP.m   &  10 &   -64.5751    &  8.5e-13   \\          \hline
Sc50b 
         &  mehrotra.m   &   8 & -70.0000    & 9.1e-07  \\
         &    arcLP.m   &  10 &   -70.0000    &  3.6e-12  \\
          \hline
Scagr25 
         &  mehrotra.m   &  18 & -1.4753e+07 & 4.6e-09  \\
         &    arcLP.m   & 19  &  -1.4753e+07  &   1.7e-08 \\ \hline
Scagr7 
         &  mehrotra.m   &  17 & -2.3314e+06 & 1.1e-07  \\
         &    arcLP.m   & 17  &   -2.3314e+06   &   7.0e-10 \\ \hline
Scfxm1+ 
         &  mehrotra.m   & 22 &  1.8417e+04 & 1.6e-08  \\
         &    arcLP.m   & 21 &  1.8417e+04 & 3.3e-05 \\ \hline
Scfxm2
        &    arcLP.m   &  24  &   3.6660e+04   &  4.8e-05   \\
         &  mehrotra.m   &   26 &  3.6660e+04 & 2.6e-08  \\ \hline
Scfxm3+
         &  mehrotra.m   &  23 &  5.4901e+04  & 9.8e-08  \\
         &    arcLP.m   &  23 &  5.4901e+04  & 1.2e-04 \\ \hline
Scrs8  
         &  mehrotra.m   & 30  &  9.0430e+02 & 1.8e-10  \\
         &    arcLP.m   & 28  &   9.0430e+02  &  1.0e-10   \\ \hline
Scsd1 
         &  mehrotra.m   &  13 &     8.6666  & 8.7e-14  \\
         &    arcLP.m   & 11  &      8.6666  & 3.3e-15  \\ \hline
Scsd6 
         &  mehrotra.m   &  16 &     50.5000 & 8.6e-15  \\
         &    arcLP.m   & 16  &   50.5000 &   2.6e-13 \\  \hline
Scsd8 
         &  mehrotra.m   &  14 &  9.0500e+02    & 1.3e-10  \\
         &    arcLP.m   & 15  &  9.0500e+02    &  2.6e-13 \\ \hline
Sctap1
         &  mehrotra.m   &  27 &  1.4123e+03   & 0.0031  \\
         &    arcLP.m   & 20  &  1.4123e+03   &   1.4e-11  \\ \hline
Sctap2 
         &  mehrotra.m   &   21 &  1.7248e+03 & 4.4e-07  \\
         &    arcLP.m   &  22  &  1.7248e+03     &  1.4e-12 \\ \hline
Sctap3
         &  mehrotra.m   &   22 &  1.4240e+03 & 5.9e-07  \\
         &    arcLP.m   &  21  &   1.4240e+03  &   1.9e-12  \\ \hline
Share1b
         &  mehrotra.m   &   25 & -7.6589e+04 & 1.5e-06  \\
         &    arcLP.m   &  26  &  -7.6589e+04  &    1.9e-07 \\ \hline
Share2b
         &  mehrotra.m   &   15 & -4.1573e+02 & 7.9e-10  \\
         &    arcLP.m   &  15  & -4.1573e+02   &   1.4e-10  \\ \hline
Ship04l
         &  mehrotra.m   &   18 &  1.7933e+06 & 2.9e-11  \\
         &    arcLP.m   &  19  &  1.7933e+06  &  1.3e-10 \\ \hline
Ship04s
         &  mehrotra.m   &   20 &  1.7987e+06 & 4.5e-09  \\
         &    arcLP.m   &  19  &    1.7987e+06   & 3.1e-10 \\ \hline
Ship08l
         &  mehrotra.m   &   22 &  1.9091e+06 & 1.0e-10  \\
         &    arcLP.m   &  20  &   1.9090e+06    &  1.8e-11  \\ \hline
Ship08s
         &  mehrotra.m   &   20 &  1.9201e+06 & 4.5e-12  \\
         &    arcLP.m   & 19   &    1.9201e+06   &   1.7e-09 \\ \hline
Ship12l
         &  mehrotra.m   &   21 &  1.4702e+06 & 1.0e-08  \\
         &    arcLP.m   & 21   &  1.4702e+06  &  3.0e-10     \\ \hline
Ship12s 
         &  mehrotra.m   &  19  &  1.4892e+06 & 2.1e-13  \\
         &    arcLP.m   & 21   &    1.4892e+06    & 5.0e-11 \\ \hline
Stocfor1+
         &  mehrotra.m  &   14  & -4.1132e+04 & 1.1e-10  \\
         &    arcLP.m   &  13  &   -4.1132e+04  & 8.6890e-11 \\ \hline
Stocfor2
         &  mehrotra.m   &   22 & -3.9024e+04 & 1.6e-09  \\
         &    arcLP.m   & 22   &  -3.9024e+04    & 4.3e-09 \\ \hline
Stocfor3
         &  mehrotra.m   & 38  & -3.9976e+04  & 6.4e-08  \\
         &    arcLP.m   & 37  & -3.9977e+04  & 7.7e-08  \\ \hline
Truss  
         &  mehrotra.m   & 26  &  4.5882e+05 & 9.5e-06  \\
         &    arcLP.m   & 24  &   4.5882e+05  & 5.2e-07   \\  \hline
\end{longtable}
\end{center}

We have an option of handling degenerate
solutions implemented in {\tt arcLP.m} and {\tt mehrotra.m}. 
For problems marked with '+', this option has to be called only by Mehrotra's method. 
For problems marked with '*', this option has to be called by both arc-search and Mehrotra's methods. 
A performance profile figure\footnote{Although the performance profile is sometimes named after Dolan-More, to our best knowledge, it was first introduced in \cite{ty96}.} is provided in 
\cite{yang18}. Simply speaking, the performance
of {\tt arcLP} is more efficient and robust than 
{\tt mehrotra}. Therefore, this solves a long-standing
dilemma in interior-point method \cite{todd02}: the best algorithm 
in theory (the short-step algorithm with the lowest/best
polynomial bound) performs poorly in numerical test, 
while the most efficient interior-point algorithm (Mehrotra's
predictor-corrector algorithm) does not have any convergence result. The 
implementation of {\tt arcLP} demonstrates 
that a theoretically attractive algorithm proposed in \cite{yang18} with
the lowest polynomial bound can also perform well 
in terms of efficiency and robustness.
This code will provide users with a free software package
whose algorithm is at least competitive to the algorithms implemented in the 
commercial software. Given broad applications of linear
programming, many users would be beneficial from
this code.


\section*{(2) Availability}
\vspace{0.5cm}
\section*{Operating system}
arcLP can be used in Linux and Windows operating systems. It has been tested in both operating systems. (A reviewer tested it on macOS operating systems.)

\section*{Programming language}
arcLP is written in Matlab. Users must have Matlab installed in their computers. Since Mathworks has compilers that can convert Matlab code to C/C++ or Fortran codes, users may take advantage of these tools to create compiled code
which can reduce the computational time by at least one order of magnitude.

\section*{Additional system requirements}
No special requirements are needed.

\section*{Dependencies}
No additional dependencies are required.

\section*{List of contributors}
Yaguang Yang is the sole author of the software.

\section*{Software location:}

{\bf Archive} GitHub.

\begin{description}[noitemsep,topsep=0pt]
	\item[Name:] arcLP
	\item[Persistent identifier:] https://github.com/yaguangyang/arcLP.git/
	\item[Licence:] \textcolor{blue}{BSD 3-clause "New" or "Revised" license.}
	\item[Publisher:]  \textcolor{blue}{Yaguang Yang.}
	\item[Version published:] \textcolor{blue}{1.0.}
	\item[Date published:] \textcolor{blue}{08/02/2026}
\end{description}

%
%
%

\section*{Language}
English

\section*{(3) Reuse potential}
Assume that the user has a Matlab installed in their Linux or Windows systems. For any linear programming problem represented in standard form with matrix A, vectors b and c, the user just needs to input or load A, b, and c in the command line and call [x,obj,kk,infe,lambda,s,exflag]=arcLP(A,b,c,d,tol,iter). This is straightforward and there is no need for any modification or extension.

\section*{Acknowledgements}
The author thanks JORS very much for kindly granting the waiver of the publication cost.

\section*{Funding statement}
JORS has generously granted the waiver of the publication fee. 

\section*{Competing interests}
The authors declare that they have no competing interests.


\vspace{2cm}

\rule{\textwidth}{1pt}

{ \bf Copyright Notice} \\
Authors who publish with this journal agree to the following terms: \\

Authors retain copyright and grant the journal right of first publication with the work simultaneously licensed under a  \href{http://creativecommons.org/licenses/by/3.0/}{Creative Commons Attribution License} that allows others to share the work with an acknowledgement of the work's authorship and initial publication in this journal. \\

Authors are able to enter into separate, additional contractual arrangements for the non-exclusive distribution of the journal's published version of the work (e.g., post it to an institutional repository or publish it in a book), with an acknowledgement of its initial publication in this journal. \\

By submitting this paper you agree to the terms of this Copyright Notice, which will apply to this submission if and when it is published by this journal.

\end{document}